\documentclass[preprint,12pt]{elsarticle}
\usepackage{amsthm}
\usepackage{amsmath}
\usepackage{amssymb}
\usepackage{amstext}
\usepackage{bm}
\usepackage{graphicx}
\usepackage{epstopdf}
\usepackage{subcaption}

 \newcommand{\bl}{\big<}
  \newcommand{\bg}{\big>}
 
\biboptions{sort&compress}

\journal{ArXiv.org}

\begin{document}

\begin{frontmatter}

\title{Adjusted Levermore-Pomraning equations for
diffusive random systems in slab geometry}

\author[rwth]{Richard Vasques\corref{cor1}}
\author[rwth]{Nitin K. Yadav\fnref{presad}}
\cortext[cor1]{Corresponding author: richard.vasques@fulbrightmail.org; Tel: +55 51 3307 9034\\
Present address: PROMEC - Federal University of Rio Grande do Sul, Brazil\\
Postal address: PROMEC-UFRGS, Av. Osvaldo Aranha 99, 90046-900 Porto Alegre, RS, Brazil}
\fntext[presad]{n.k.yadav@tue.nl; Tel: +31 617 38 9780\\
Present address: Department of Mathematics and Computer Science; Centre for Analysis, Scientific Computing and Applications; Eindhoven University of Technology, Netherlands}
\address[rwth]{Department of Mathematics,
Center for Computational Engineering Science,\\
RWTH Aachen University, Aachen, Germany}

\begin{abstract}
This paper presents a multiple length-scale asymptotic analysis
for transport problems in 1-D diffusive random media. This analysis shows that the Levermore-Pomraning (LP) equations can be adjusted in order to achieve the correct asymptotic behavior. This adjustment appears in the form of a rescaling of the Markov transition functions by a factor $\eta$, which can be chosen in a simple way. Numerical results are given that (i) validate the theoretical predictions; and (ii) show that the adjusted LP equations greatly outperform the standard LP model for this class of transport problems. 

\end{abstract}


\end{frontmatter}

\section{Introduction}\label{sec1}
\setcounter{section}{1}
\setcounter{equation}{0} 

The {\em diffusion approximation} is a classic model for particle transport in a physical system in which absorption and sources are weak and the solution varies slowly over the distance of a mean
free path. The diffusion equation has been shown to be an asymptotic limit of the transport equation
\cite{larsen74,habetler,larsen80}.

The {\em Levermore-Pomraning (LP) equations}
\cite{adams, pomraning91, pomraning98a, pomraning98b}
are a well-known approach to model particle transport in a heterogeneous physical system consisting of two or more materials. They generalize the widely-used {\em atomic mix model}
\cite{pomraning91,dumas},
which requires chunks of each material to be optically thin. The LP model is known to be accurate
for problems with (i) weak scattering, and (ii) a Markovian distribution of chord lengths across material
chunks \cite{vasques04,brantley}.
However, numerical experiments have indicated that the LP model is inaccurate for diffusive systems \cite{vasques04,davis}.

In this work, we consider 1-D transport problems in slab geometry. In particular, we assume that:
\begin{itemize}
\item[{\bf A1.}] The physical system is heterogeneous, consisting of alternate layers of different materials. The thicknesses of the layers are on the order
of a mean free path (or smaller).
\item[{\bf A2.}]For convenience, we assume that the system is binary, with different layers labeled 1 and 2. The cross sections and source for material $i$ ($i$ = 1 or 2) are labeled $\Sigma_{ti}$, $\Sigma_{si}$, and $Q_i$.
\item[{\bf A3.}]The geometrical structure of the system is a stochastic binary mixture, with the mixing statistics taken as Markovian. 
\item[{\bf A4.}] The system is {\em diffusive} (in a global sense). That is, the physical system is optically thick, and absorption and sources are weak at each spatial point.
\end{itemize}
We point out that assumptions {\bf A1} and {\bf A4} imply that the number of material layers in the system is large.

As introduced in \cite{larsen05}, for the diffusive limit just described, the standard LP model reduces to a diffusion
equation with an incorrect diffusion coefficient. This theoretically explains the inaccuracies observed in LP simulations for diffusive problems. In this paper, we present an asymptotic analysis that leads to an {\em Adjusted LP model} (ALP); this model preserves the correct asymptotic behavior in this diffusive limit, greatly improving the accuracy of the LP equations.

We include in this work numerical simulations that validate the asymptotic theory presented for this diffusive limit. These simulations demonstrate that: (i) the solution of the LP equations limits to the solution of the LP diffusion equation predicted by our asymptotic theory; (ii) the solution of the ALP equations limits to the solution of the ALP diffusion equation predicted by our asymptotic theory; (iii) the solution of the ALP model accurately agrees with the correct (ensemble-averaged) solution of the transport problem, while the standard LP model fails to do the same.

We also include results showing that the ALP model preserves the accuracy of the standard LP equations for problems with weak or no scattering. Furthermore, we present numerical simulations indicating that the ALP equations remain accurate for problems in which the diffusive characteristics of the medium are slightly relaxed. In these simulations, the ALP model is shown to outperform both its standard counterpart {\em and} the atomic mix model.  

A summary of the remainder of the paper follows.
In Section \ref{sec2} we present the asymptotic theory for the
transport equation in the 1-D random diffusive system. In Section \ref{sec3} we present the LP model and propose the Adjusted LP equations, introducing a factor $\eta$ to be defined later. We perform an asymptotic analysis of the ALP equations for $\eta$ of $O(1)$ (Section \ref{sec3.1}) and for $\eta$ of $O(1/\varepsilon)$ (Section \ref{sec3.2}).
In Section \ref{sec4} we propose an expression for $\eta$ motivated by the analysis performed in Section \ref{sec3}. We
present numerical results for problems in diffusive (Section \ref{sec4.1}) and non-diffusive (Section \ref{sec4.2}) systems. These results confirm the predictions of the asymptotic analysis and validate the ALP model. We conclude with a discussion in Section \ref{sec5}.

\section{Asymptotic Analysis of the Transport Equation}\label{sec2}
\setcounter{section}{2}
\setcounter{equation}{0} 

\begin{figure}
    \centering
        \includegraphics[scale=0.4]{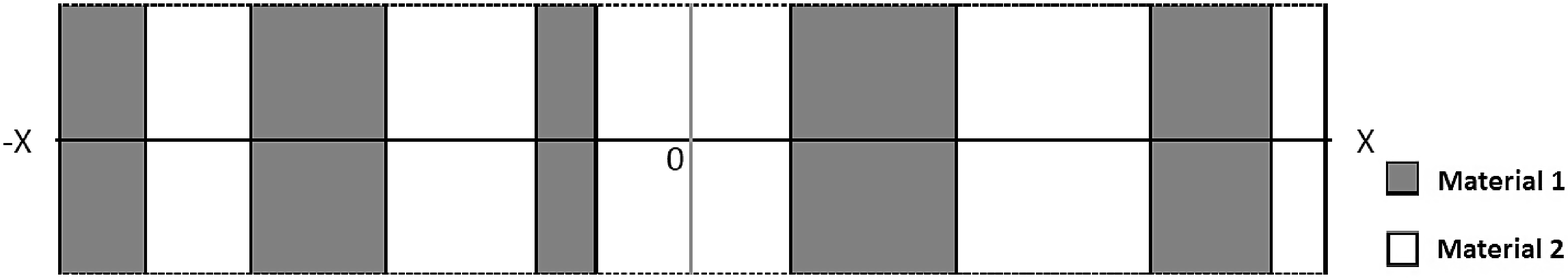}
            \caption{Sketch of the 1-D stochastic system}\label{fig1}
\end{figure}

 We consider the following 1-D steady-state, monoenergetic transport problem, with vacuum boundaries and space-dependent cross sections and source:
\begin{subequations}\label{2.1}
\begin{align}
\mu\frac{\partial\psi}{\partial x}(x,\mu) + \Sigma_t(x)\psi(x,\mu) &=\frac{\Sigma_s(x)}{2}\int_{-1}^1\psi(x,\mu')d\mu'+\frac{Q(x)}{2}, \label{2.1a}
\\& \hspace{4cm} -X\leq x\leq X, -1\leq\mu\leq 1, \nonumber
\\& \psi(-X,\mu)=0, \,\,\, 0<\mu\leq 1,
\\& \psi(X,\mu)=0, \,\,\, -1\leq\mu< 0.
\end{align}
\end{subequations}
Moreover, taking into account the assumptions described in the previous section, we consider the following:
\begin{itemize}
\item[I.] The physical system $-X\leq x\leq X$ consists of a stochastic structure of alternate layers of two distinct materials; a realization of the system is sketched in Figure \ref{fig1}. The cross sections and source in Eq.\ (\ref{2.1a}) are stochastic functions of space.
\item[II.] The spatial variable $x$ is scaled so that the typical width of a layer and a typical mean free path are $O(1)$. Thus, $\Sigma_t = O(1)$.
\item[III.] The system is optically thick. Thus, the dimensionless parameter
\begin{subequations}\label{2.2}
\begin{align}
\varepsilon \equiv \frac{\text{average width of a layer}}{\text{width of the system}} = \frac{1}{\text{number of layers}}
\end{align}
is small, and
\begin{align}
2X = \text{width of the system} = O(1/\varepsilon).
\end{align}
\end{subequations}
\item[IV.] Absorption is weak at all spatial points. This is expressed by writing the absorption cross section as
\begin{equation}\label{2.3}
\Sigma_t(x)-\Sigma_s(x) = \Sigma_a(x) = \varepsilon^2\sigma_a(x),
\end{equation}
where $\sigma_a(x) = O(1)$.
\item[V.] For convenience, we scale the source so that the infinite-medium solution is $O(1)$. This is expressed by writing
\begin{equation}\label{2.4}
Q(x) = \varepsilon^2 q(x),
\end{equation}
where $q(x) = O(1)$.
\item[VI.] The flux depends on two spatial variables: the ``fast" spatial variable $x$, which describes rapid variations on the order of a mean free path or a layer width, and a new ``slow" spatial variable
\begin{equation}\label{2.5}
z=\varepsilon x,
\end{equation}
which describes slowly-varying spatial variations in the flux over the $O(1/\varepsilon)$ optical width of the system. This is expressed by writing
\begin{subequations}\label{2.6}
\begin{align}
\psi(x,\mu) = \hat\psi(x,z,\mu),
\end{align}
which implies:
\begin{align}
\frac{\partial\psi}{\partial x}(x,\mu) = \frac{\partial\hat\psi}{\partial x}(x,z,\mu) +\varepsilon \frac{\partial\hat\psi}{\partial z}(x,z,\mu).
\end{align}
\end{subequations}
\end{itemize}
Introducing Eqs.\ (\ref{2.3})-(\ref{2.6}) into Eq.\ (\ref{2.1a}), we obtain the following scaled transport equation:
\begin{equation}\label{2.7}
\begin{split}
\mu\frac{\partial \hat\psi}{\partial x}(x,z,\mu)+\varepsilon\mu\frac{\partial\hat\psi}{\partial z}(x,z,\mu)+&\Sigma_t(x)\hat\psi(x,z,\mu)
\\& =\frac{\Sigma_t(x)-\varepsilon^2\sigma_a(x)}{2}\int_{-1}^1\hat\psi(x,z,\mu')d\mu'+\varepsilon^2\frac{q(x)}{2}.
\end{split}
\end{equation}
This equation can be asymptotically solved by assuming the ansatz
 \begin{equation}
 \hat\psi(x,z,\mu) = \sum_{n=0}^\infty \varepsilon^n\hat\psi_n(x,z,\mu),
 \end{equation}   
 in which $\varepsilon << 1$, and $x$ and $z$ are treated as independent variables. This calculation, given in detail elsewhere
\cite{larsen05},
yields the diffusion equation 
 \begin{align}\label{2.9}
-\frac{1}{3\bl\Sigma_t\bg}\frac{d^2\phi_0}{dz^2}(z)+\bl\sigma_a\bg\phi_0(z)=\bl q\bg,
\end{align}
where
\begin{subequations}
\begin{align}
\bl\Sigma_t\bg &= \text{volume-averaged total cross section},
\\\bl\sigma_a\bg &= \text{volume-averaged absorption cross section},
\\\bl q\bg &= \text{volume-averaged interior source}.
\end{align}
\end{subequations} 
Finally, we return to the original unstrecthed variables. Multiplying Eq.\ (\ref{2.9}) by $\varepsilon^2$ and using
\begin{align}\label{2.11}
\Phi_0(x) &= \phi_0(\varepsilon x) = \phi_0(z),
\end{align}
we obtain
\begin{equation}\label{2.12}
-\frac{1}{3\bl\Sigma_t\bg}\frac{d^2\Phi_0}{dx^2}(x)+\bl\Sigma_a\bg\Phi_0(x)=\bl Q\bg, \,\,\,\, -X<x<X\, .
\end{equation}
Notice that Eq.\ (\ref{2.12}) is the leading-order asymptotic limit of Eq.\ (\ref{2.7}) as $\varepsilon\rightarrow 0$, and the unknown $\Phi_0(x)$ is the leading-order estimate of the scalar flux.

\textbf{Note:} The atomic mix approximation \cite{pomraning91,dumas} of Eq.\ (\ref{2.1a}) is
\begin{align}\label{2.13}
\mu\frac{\partial\psi}{\partial x}(x,\mu) + \bl\Sigma_t\bg\psi(x,\mu) &=\frac{\bl\Sigma_s\bg}{2}\int_{-1}^1\psi(x,\mu')d\mu'+\frac{\bl Q\bg}{2}, 
\\& \hspace{4cm} -X\leq x\leq X, -1\leq\mu\leq 1. \nonumber
\end{align}
If the conventional diffusion approximation were
applied to Eq.\ (\ref{2.13}), the resulting equation would be the same as Eq.\ (\ref{2.12}), which is simply the conventional diffusion equation with atomic mix (volume-averaged) cross sections.

\section{The Adjusted Levermore-Pomraning (ALP) Equations}\label{sec3}
\setcounter{section}{3}
\setcounter{equation}{0} 

The Levermore-Pomraning formulation \cite{adams, pomraning91, pomraning98a, pomraning98b} for Eq.\ (\ref{2.1a}) is given by
 \begin{align}
\mu\frac{\partial p_i\Psi_i}{\partial x}&(x,\mu) + \Sigma_{ti}p_i\Psi_i(x,\mu) = \frac{\Sigma_{si}}{2}\int_{-1}^1p_i\Psi_i(x,\mu')d\mu' \label{3.1}
\\& +\left(\frac{p_j\overline\Psi_j(x,\mu)}{\Lambda_j(x,\mu)}-\frac{p_i\overline\Psi_i(x,\mu)}{\Lambda_i(x,\mu)}\right) +\frac{p_iQ_i}{2}, \nonumber
\,\,\,\,-X\leq x\leq X, -1\leq\mu\leq 1,
\end{align}
where $i, j$ = 1 or 2 with $j\not=i$, and
\begin{itemize}
\item[I.] $\Psi_i(x,\mu)$ is the ensemble average of $\psi(x,\mu)$ over all physical realizations such that $x$ is in material $i$;
\item[II.] $\overline\Psi_i(x,\mu)$ is the ensemble average of $\psi(x,\mu)$ over all physical realizations such that $x$ is at an interface point and $\mu$ points out of material $i$;
\item[III.] $p_i$ is the probability of finding material $i$ at point $x$;
\item[IV.] $\Lambda_i(x,\mu)$ is the {\em Markov transition function} of material $i$, defined such that a particle moving a distance $ds$ in a direction $\mu$ has a probability $ds/\Lambda_i(x,\mu)$ of transferring from material $i$ to material $j\neq i$.
\end{itemize}
For Markovian statistics, $\Lambda_i(x,\mu)$ is simply the mean chord length in direction $\mu$ in material $i$. Therefore, defining
\begin{subequations}\label{3.2}
\begin{align}\label{3.2a}
\lambda_i &= \text{the mean width of the layers of material $i$}, 
\end{align}
we write
\begin{align}
&\Lambda_i(x,\mu) = \frac{\lambda_i}{|\mu|},
\\ &p_i = \frac{\lambda_i}{\lambda_1+\lambda_2} = \text{the volume fraction of material $i$} .\label{3.2c}
\end{align}
\end{subequations}
It is important to notice that Eq.\ (\ref{3.1}) represents two equations with four unknowns, namely $\Psi_1, \Psi_2, \overline\Psi_1$, and $\overline\Psi_2$. To obtain a useful set of equations, a closure is introduced: $\overline\Psi_i$ is replaced with $\Psi_i$. 

Using Eqs.\ (\ref{3.2}), the {\em standard} LP model for Eqs.\ (\ref{2.1}) is written as
\begin{subequations}\label{3.3}
 \begin{align}
\mu\frac{\partial p_i\Psi_i}{\partial x}&(x,\mu) + \Sigma_{ti}p_i\Psi_i(x,\mu) = \frac{\Sigma_{si}}{2}\int_{-1}^1p_i\Psi_i(x,\mu')d\mu' 
\\& +|\mu|\left(\frac{p_j\Psi_j(x,\mu)}{\lambda_j}-\frac{p_i\Psi_i(x,\mu)}{\lambda_i}\right) +\frac{p_iQ_i}{2}, \nonumber
\,\,\,\,-X\leq x\leq X, -1\leq\mu\leq 1,
\\& \hspace{2cm} \Psi_i(-X,\mu)=0, \,\,\, 0<\mu\leq 1,
\\& \hspace{2cm} \Psi_i(X,\mu)=0, \,\,\, -1\leq\mu< 0,
\end{align}
such that the LP estimate of $\bl \Psi \bg$ (the ensemble average of the mean angular flux over all physical realizations) is given by
\begin{align}\label{3.3d}
\bl\Psi\bg(x,\mu)=p_1\Psi_1(x,\mu) +p_2\Psi_2(x,\mu).
\end{align}
\end{subequations}
Equations (\ref{3.3}) are known to model problems in purely absorbing media exactly, and to be accurate for problems with weak scattering. However, the accuracy of this formulation decreases as scattering increases \cite{vasques04,davis}.

In Section \ref{sec2}, we have shown that 
the solution of the transport equation (\ref{2.1a}) limits to the solution of the (atomic mix) diffusion equation (\ref{2.12}) in the 1-D diffusive limit. Moreover, when $\Lambda_i\rightarrow 0$, it has been shown \cite{pomraning91} that Eqs.\ (\ref{3.3}) limit to the atomic mix formulation given by Eq.\ (\ref{2.13}). 
This raises the following question:
\begin{itemize}
\item {\em Is it possible to scale the Markov transition functions $\Lambda_i(x,\mu)$ in such a way that: (i) the resulting Adjusted LP equations yield the correct diffusion equation given in Eq.\ (\ref{2.12}); and (ii) the standard LP equations are preserved for purely absorbing media?}
\end{itemize}
To answer this question, let us scale the Markov transition functions such that
\begin{align}\label{3.4}
\Lambda_i(x,\mu) = \frac{1}{\eta}\frac{\lambda_i}{|\mu|},
\end{align}
and let us write the ALP equations as
\begin{subequations}\label{3.5}
 \begin{align}\label{3.5a}
\mu\frac{\partial p_i\Psi_i}{\partial x}&(x,\mu) + \Sigma_{ti}p_i\Psi_i(x,\mu) = \frac{\Sigma_{si}}{2}\int_{-1}^1p_i\Psi_i(x,\mu')d\mu' 
\\& +{\eta}|\mu|\left(\frac{p_j\Psi_j(x,\mu)}{\lambda_j}-\frac{p_i\Psi_i(x,\mu)}{\lambda_i}\right) +\frac{p_iQ_i}{2}, \nonumber
\,\,\,\,-X\leq x\leq X, -1\leq\mu\leq 1,
\\& \hspace{2cm} \Psi_i(-X,\mu)=0, \,\,\, 0<\mu\leq 1,
\\& \hspace{2cm} \Psi_i(X,\mu)=0, \,\,\, -1\leq\mu< 0,
\end{align}
\end{subequations}
where the factor $1/\eta$ was taken outside of the parenthesis. 
Our aim is to find an expression for $\eta$ that (i) yields accurate results in diffusive media; and (ii) yields $\eta=1$ in purely absorbing media, in which case Eqs.\ (\ref{3.5}) reduce to the correct standard LP model given by Eqs.\ (\ref{3.3}).

\subsection{Asymptotic Analysis of the ALP Equations with $\eta$ of $O(1)$}\label{sec3.1}

To analyze Eqs.\ (\ref{3.5}) in the same asymptotic limit applied to Eqs.\ (\ref{2.1}), we simply take $\Sigma_{ti}$ and $\lambda_i$ to be $O(1)$, $X$ to be $O(1/\varepsilon)$, and:
\begin{subequations}\label{3.6}
 \begin{align}
\Sigma_{ai}&=\Sigma_{ti}-\Sigma_{si} =\varepsilon^2\sigma_{ai},
\\ Q_i &=\varepsilon^2q_i ,
\\ \Psi_i(x,\mu)&=\psi_i(z,\mu), \label{3.6c}
\end{align}
\end{subequations}
where $z$ is given by Eq.\ (\ref{2.5}).
Here, the flux depends only on the ``slow" spatial variable $z$ and the angular variable $\mu$, there being no ``fast" spatial variation in Eqs.\ (\ref{3.5}).

Assuming $\eta$ to be $O(1)$ and introducing Eqs.\ (\ref{3.6}) into Eq.\ (\ref{3.5a}), we obtain the scaled ALP equations:
\begin{equation}\label{3.7}
\begin{split}
\varepsilon\mu\frac{\partial p_i\psi_i}{\partial z}(z,\mu) + \Sigma_{ti}p_i\psi_i&(z,\mu) = 
\frac{\Sigma_{ti}-\varepsilon^2\sigma_{ai}}{2}\int_{-1}^1p_i\psi_i(z,\mu')d\mu'
\\&+ \eta|\mu|\left(\frac{p_j\psi_j(z,\mu)}{\lambda_j}-\frac{p_i\psi_i(z,\mu)}{\lambda_i}\right) +\varepsilon^2\frac{p_iq_i}{2}.
\end{split}
\end{equation}
We solve these equations by assuming the ansatz
\begin{equation}\label{3.8}
 \psi_i(z,\mu) = \sum_{n=0}^\infty \varepsilon^n\psi_{i,n}(z,\mu).
 \end{equation}   
Introducing Eqs.\ (\ref{3.8}) into Eq.\ (\ref{3.7}) and equating the coefficients of different powers of $\varepsilon$, we obtain for $n \geq 0$:
\begin{equation}\label{3.9}
\begin{split}
\Sigma_{ti}&\left[p_i\psi_{i,n}(z,\mu)- \frac{1}{2}\int_{-1}^1p_i\psi_{i,n}(z,\mu')d\mu'\right]+ \eta|\mu|\left(\frac{p_i\psi_{i,n}(z,\mu)}{\lambda_i}-\frac{p_j\psi_{j,n}(z,\mu)}{\lambda_j}\right)
\\&\hspace{2.8cm}=-\mu\frac{\partial p_i\psi_{i,n-1}}{\partial z}(z,\mu)
-\frac{\sigma_{ai}}{2}\int_{-1}^1p_i\psi_{i,n-2}(z,\mu')d\mu'+\delta_{n,2}\frac{p_iq_i}{2},
\end{split}
\end{equation}
￼￼￼where $\psi_{i,-1}=\psi_{i,-2}=0$. These equations can be solved recursively:
for $n = 0$, they have only an isotropic solution of the form
\begin{equation}\label{3.10}
\psi_{i,0}(z,\mu) = \frac{\phi_0(z)}{2},
\end{equation}
where $\phi_0(z)$ is undetermined. 

For $n = 1$, Eqs.\ (\ref{3.9}) have a solvability condition that is automatically satisfied. The general solution of the $n = 1$ equations is:
\begin{equation}\label{3.11}
\psi_{i,1}(z,\mu) = \frac{1}{2}\left[\phi_1(z)-\mu f_i(|\mu|)\frac{d\phi_0}{dz}(z)\right],
\end{equation}
￼￼where $\phi_1(z)$ is undetermined and:
\begin{equation}
f_i(|\mu|) = \frac{\lambda_1\lambda_2\Sigma_{tj}+(\lambda_1+\lambda_2)\eta|\mu|}{\lambda_1\lambda_2\Sigma_{t1}\Sigma_{t2} + (\lambda_1\Sigma_{t1}+\lambda_2\Sigma_{t2})\eta|\mu|}.
\end{equation}

For $n = 2$,  Eqs.\ (\ref{3.9}) have a solvability condition that is {\em not} automatically satisfied. This condition is obtained by first integrating Eqs.\ (\ref{3.9}) with $n = 2$ over $-1 \leq\mu\leq 1$, and then adding the resulting two equations.
This gives:
\begin{equation}\label{3.13}
\begin{split}
0 = -\frac{d}{dz}\int_{-1}^1&\mu[p_1\psi_{1,1}(z,\mu)+p_2\psi_{2,1}(z,\mu)]d\mu
\\& -\int_{-1}^1 [p_1\sigma_{a1}\psi_{1,0}(z,\mu)+p_2\sigma_{a2}\psi_{2,0}(z,\mu)]d\mu
 + (p_1q_1+p_2q_2).
\end{split}
\end{equation}
Introducing Eqs.\ (\ref{3.10}) and (\ref{3.11}) into Eq.\ (\ref{3.13}) and simplifying, we obtain the following diffusion equation for $\phi_0$:
\begin{equation}\label{3.14}
-\frac{\beta}{3\bl\Sigma_t\bg}\frac{d^2\phi_0}{dz^2}(z)+\bl\sigma_a\bg\phi_0(z)=\bl q\bg,
\end{equation}
where
\begin{subequations}
\begin{align}
\bl\sigma_a\bg &= p_1\sigma_{a1} + p_2\sigma_{a2} ,
\\ \bl q\bg &= p_1q_1 + p_2q_2 ,
\\ \beta &= \int_0^1 3\mu^2 \alpha(\mu)d\mu, \label{3.15c}
\\ \alpha(\mu) &=\frac{\lambda_1\lambda_2\bl\Sigma_t\bg(p_1\Sigma_{t2}+p_2\Sigma_{t1})
+\eta(\lambda_1\Sigma_{t1}+\lambda_2\Sigma_{t2})\mu}{\lambda_1\lambda_2\Sigma_{t1}\Sigma_{t2} + \eta(\lambda_1\Sigma_{t1}+\lambda_2\Sigma_{t2})\mu}.
\end{align}
\end{subequations}
Multiplying Eq.\ (\ref{3.14}) by $\varepsilon^2$ we obtain:
\begin{equation}\label{3.16}
-\frac{\beta}{3\bl\Sigma_t\bg}\frac{d^2\Phi_0}{dx^2}(x)+\bl\Sigma_a\bg\Phi_0(x)=\bl Q\bg.
\end{equation}
Finally, Eqs.\ (\ref{2.11}), (\ref{3.2c}), (\ref{3.3d}), (\ref{3.6c}), (\ref{3.8}), and (\ref{3.10}) give:
\begin{equation}
\Phi(x) = \int_{-1}^1\bl\Psi\bg(x,\mu)d\mu = \Phi_0(x) + O(\varepsilon).
\end{equation}
￼￼￼￼￼Thus, Eq.\ (\ref{3.16}) is the leading-order asymptotic limit of Eq.\ (\ref{3.7}) as $\varepsilon \rightarrow 0$, and the unknown $\Phi_0(x)$ in this equation is the leading-order estimate of the scalar flux.

We note that the diffusion coefficient in Eq.\ (\ref{3.16}) is equal to the one in Eq.\ (\ref{2.12}) {\em only if} $\beta=1$, which only happens if $\alpha(\mu)=1$. However, if $\Sigma_{t1}\neq\Sigma_{t2}$, it is easy to show that $\alpha(\mu) > 1$ for any choice of $\lambda_i$, which implies that $\beta > 1$. This causes the diffusion coefficient of Eq.\ (\ref{3.16}) to be unphysically large, which considerably affects the solution. Numerical results presented in Section \ref{sec4} show that this is the case for the standard LP equations, in which $\eta=1$.

Furthermore, this analysis indicates that if we choose $\eta>>1$ such that $\Lambda_i<< 1$, both $\alpha(\mu)$ and $\beta$ approach $1$. This suggests that $\eta$ should be chosen as $O(1/\varepsilon^k)$ for some $k>0$. In the next section, we present an asymptotic analysis for Eqs.\ (\ref{3.5}) with $\eta$ chosen as $O(1/\varepsilon)$.    

\subsection{Asymptotic Analysis of the ALP Equations with $\eta$ of $O(1/\varepsilon)$}\label{sec3.2}

Let us define the change of variables
\begin{subequations}
\begin{align}
\bl\Psi\bg(x,\mu) &= p_1\Psi_1(x,\mu) + p_2\Psi_2(x,\mu),\\
\vartheta(x,\mu)&= \sqrt{p_1p_1}\,[\Psi_1(x,\mu) - \Psi_2(x,\mu)],
\end{align}
\end{subequations}
and rewrite Eqs.\ (\ref{3.5a}) in an algebraically different (but equivalent) form \cite{pomraning91}:
\begin{align}\label{3.19}
\mu\frac{\partial}{\partial x}
\left[\begin{array}{c}\bl\Psi\bg\\ \vartheta\end{array}\right]
+
\left[\begin{array}{cc}\bl\Sigma_t\bg & \nu \\ \nu & \hat\Sigma_t\end{array}\right]
\left[\begin{array}{c}\bl\Psi\bg\\ \vartheta\end{array}\right]
=
\frac{1}{2}
\left[\begin{array}{cc}\bl\Sigma_s\bg & \nu_s \\ \nu_s & \hat\Sigma_s\end{array}\right]
\left[\begin{array}{c}\Phi\\ \varphi\end{array}\right]
+
\frac{1}{2}\left[\begin{array}{c}\bl Q\bg\\ U\end{array}\right].
\end{align}
Here, $\bl\Sigma_t\bg$, $\bl\Sigma_s\bg$, and $\bl Q\bg$ are the volume-averaged cross sections and source, and
\begin{subequations}
\begin{align}
\hat\Sigma_t &= \hat\Sigma_t(|\mu|) =  p_2\Sigma_{t1} + p_1\Sigma_{t2} + \eta |\mu| \frac{\lambda_1+\lambda_2}{\lambda_1\lambda_2}, \label{3.20a}\\
\hat\Sigma_s &= p_2\Sigma_{s1} + p_1\Sigma_{s2},\\
\nu &= \sqrt{p_1p_2}\,(\Sigma_{t1}-\Sigma_{t2}),\\
\nu_s &= \sqrt{p_1p_2}\,(\Sigma_{s1}-\Sigma_{s2}),\\
U &=\sqrt{p_1p_2}\,(Q_1-Q_2),\\
\Phi &= \Phi(x) = \int_{-1}^1\bl\Psi\bg(x,\mu)d\mu,\\
\varphi &= \varphi(x) = \int_{-1}^1\vartheta(x,\mu)d\mu.
\end{align}
\end{subequations}

As we have done in Section \ref{sec3.1}, we take $\Sigma_{ti}$ and $\lambda_i$ to be $O(1)$, $X$ to be $O(1/\varepsilon)$, and define
\begin{subequations}\label{3.21}
 \begin{align}
\bl\Sigma_s\bg &= \bl\Sigma_t\bg-\bl\Sigma_a\bg=\bl\Sigma_t\bg-\varepsilon^2\bl\sigma_a\bg,\\
\nu_s &= \sqrt{p_1p_2}\, [(\Sigma_{t1}-\Sigma_{t2})-(\varepsilon^2\sigma_{a1}-\varepsilon^2\sigma_{a2})] = \nu-\varepsilon^2\nu_a,\\
\bl Q\bg &=\varepsilon^2\bl q\bg ,\\
U &=\sqrt{p_1p_2}\, (\varepsilon^2q_1-\varepsilon^2q_2) = \varepsilon^2 u,\\
\bl \Psi\bg(x,\mu)&=\bl \Psi\bg(\varepsilon z,\mu)=\psi(z,\mu),\\
\vartheta(x,\mu)&=\vartheta(\varepsilon z,\mu)=\overline\vartheta(z,\mu),\\
\Phi (x) &= \Phi (\varepsilon z) = \phi(z),\\
\varphi(x) &= \varphi(\varepsilon z) = \overline\varphi(z).
\end{align}
Finally, we assume $\eta$ to be $O(1/\varepsilon)$. This means that $\hat\Sigma_t$ given in Eq.\ (\ref{3.20a}) is $O(1/\varepsilon)$, and therefore we define
\begin{align}
\hat\Sigma_t &= \varepsilon^{-1}\hat\sigma_t.
\end{align} 
\end{subequations}
Introducing Eqs.\ (\ref{3.21}) into Eq.\ (\ref{3.19}), we obtain
\begin{align}\label{3.22}
\varepsilon\mu\frac{\partial}{\partial z}
\left[\begin{array}{c}\psi\\ \overline\vartheta\end{array}\right]
&+
\left[\begin{array}{cc}\bl\Sigma_t\bg & \nu \\ \nu & \hat\sigma_t/\varepsilon\end{array}\right]
\left[\begin{array}{c}\psi\\ \overline\vartheta\end{array}\right]\\
&\hspace{30pt}=
\frac{1}{2}
\left[\begin{array}{cc}\bl\Sigma_t\bg-\varepsilon^2\bl\sigma_a\bg & \nu-\varepsilon^2\nu_a \\ \nu-\varepsilon^2\nu_a & \hat\Sigma_s\end{array}\right]
\left[\begin{array}{c}\phi\\ \overline\varphi\end{array}\right]
+\frac{\varepsilon^2}{2}
\left[\begin{array}{c}\bl q\bg\\ u\end{array}\right].\nonumber
\end{align}
In order to solve this equation, we introduce the asymptotic expansions
 \begin{subequations}
 \begin{align}
 \psi(z,\mu) &= \sum_{n=0}^\infty \varepsilon^n\psi_n(z,\mu),\\
 \overline\vartheta(z,\mu) &= \sum_{n=0}^\infty \varepsilon^n\overline\vartheta_n(z,\mu),
 \end{align}
 \end{subequations}   
into Eq.\ (\ref{3.22}) and equate the coefficients of different powers of $\varepsilon$.
For terms of $O(\varepsilon^{-1})$, we easily obtain
\begin{align}\label{3.24}
\overline\vartheta_0(z,\mu) = 0.
\end{align}

For terms of $O(\varepsilon^0)$, we have 
\begin{subequations}\label{3.25}
\begin{align}
\bl\Sigma_t\bg\left[\psi_0(z,\mu)-\frac{\phi_0(z)}{2}\right] + \nu\left[\overline\vartheta_0(z,\mu) -\frac{\overline\varphi_0(z)}{2}\right]= 0,\\
\nu\left[\psi_0(z,\mu)-\frac{\phi_0(z)}{2}\right] + \hat\sigma_t\overline\vartheta_1(z,\mu) = \frac{\hat\Sigma_s}{2}\overline\varphi_0(z).
\end{align}
\end{subequations}
Bearing in mind that $\overline\varphi_n(z) = \int_{-1}^1\overline\vartheta_n(z,\mu)d\mu$, we introduce Eq.\ (\ref{3.24}) into Eqs.\ (\ref{3.25}) and obtain
\begin{subequations}\label{3.26}
\begin{align}
\psi_0(z,\mu) &= \frac{\phi_0(z)}{2},\\
\overline\vartheta_1(z,\mu) &= 0,
\end{align}
\end{subequations}
￼￼where $\phi_0(z)$ is undetermined.

For terms of $O(\varepsilon^1)$, we have
\begin{subequations}\label{3.27}
\begin{align}
\bl\Sigma_t\bg\left[\psi_1(z,\mu)-\frac{\phi_1(z)}{2}\right] + \nu\left[\overline\vartheta_1(z,\mu) -\frac{\overline\varphi_1(z)}{2}\right]= -\mu\frac{\partial\psi_0}{\partial z}(z,\mu),\label{3.27a}\\
\nu\left[\psi_1(z,\mu)-\frac{\phi_1(z)}{2}\right] + \hat\sigma_t\overline\vartheta_2(z,\mu) = \frac{\hat\Sigma_s}{2}\overline\varphi_1(z) -\mu\frac{\partial\overline\vartheta_0}{\partial z}(z,\mu).\label{3.27b}
\end{align}
\end{subequations}
These equations have a solvability condition that requires $\int_{-1}^1\hat\sigma_t\overline\vartheta_2d\mu = 0$. Introducing Eqs.\ (\ref{3.24}) and (\ref{3.26}) into Eqs.\ (\ref{3.27}), we obtain
\begin{subequations}\label{3.28}
\begin{align}
\psi_1(z,\mu) &= \frac{\phi_1(z)}{2} - \frac{\mu}{2\bl\Sigma_t\bg}\frac{d\phi_0}{dz}(z),\\
\overline\vartheta_2(z,\mu) &= \frac{\mu}{2\bl\Sigma_t\bg}\frac{\nu}{\hat\sigma_t}\frac{d\phi_0}{dz}(z),
\end{align}
\end{subequations}
￼￼where $\phi_1(z)$ is undetermined. Since $\overline\vartheta_2$ is an odd function of $\mu$ and $\hat\sigma_t = \hat\sigma_t(|\mu|)$, the
solvability condition is satisfied.

For terms of $O(\varepsilon^2)$, we have
\begin{subequations}
\begin{align}
&\bl\Sigma_t\bg\left[\psi_2(z,\mu)-\frac{\phi_2(z)}{2}\right] + \nu\left[\overline\vartheta_2(z,\mu) -\frac{\overline\varphi_2(z)}{2}\right]\\
&\hspace{100pt}= -\mu\frac{\partial\psi_1}{\partial z}(z,\mu)
-\frac{\bl\sigma_a\bg}{2}\phi_0(z)-\frac{\nu_a}{2}\overline\varphi_0(z) + \frac{\bl q\bg}{2},\nonumber\\
&\nu\left[\psi_2(z,\mu)-\frac{\phi_2(z)}{2}\right] + \hat\sigma_t\overline\vartheta_3(z,\mu) \\
&\hspace{100pt}= \frac{\hat\Sigma_s}{2}\overline\varphi_2(z) -\mu\frac{\partial\overline\vartheta_1}{\partial z}(z,\mu) -\frac{\nu_a}{2}\phi_0(z)-\frac{\bl\sigma_a\bg}{2}\overline\varphi_0(z) + \frac{u}{2}\nonumber.
\end{align}
\end{subequations}
These equations have a solvability condition obtained by integrating over $-1\leq \mu\leq 1$, which gives 
\begin{subequations}\label{3.30}
\begin{align}
0&= -\frac{\partial}{\partial z}\int_{-1}^1\mu\psi_1(z,\mu)d\mu
-\bl\sigma_a\bg\phi_0(z)-\nu_a\overline\varphi_0(z) + \bl q\bg,
\\
\int_{-1}^1\hat\sigma_t\overline\vartheta_3(z,\mu)d\mu &= \hat\Sigma_s\overline\varphi_2(z) -\frac{\partial}{\partial z}\int_{-1}^1\mu\overline\vartheta_1(z,\mu)d\mu -\nu_a\phi_0(z)-\bl\sigma_a\bg\overline\varphi_0(z) + u.
\end{align}
\end{subequations}
Introducing Eqs.\ (\ref{3.24}), (\ref{3.26}), and (\ref{3.28}) into Eqs.\ (\ref{3.30}), we obtain
\begin{subequations}\label{3.31}
\begin{align}
0&= \frac{1}{3\bl\Sigma_t\bg}\frac{d^2\phi_0}{dz^2}(z)
-\bl\sigma_a\bg\phi_0(z)+ \bl q\bg,\label{3.31a}
\\
\int_{-1}^1\hat\sigma_t\overline\vartheta_3(z,\mu)d\mu &= -\nu_a\phi_0(z) + u.
\end{align}
\end{subequations}

Equation (\ref{3.31a}) is the diffusion equation for $\phi_0(z)$; multiplying this equation by $\varepsilon^2$, we obtain
\begin{align}\label{3.32}
-\frac{1}{3\bl\Sigma_t\bg}\frac{d^2\Phi_0}{dx^2}(x)
-\bl\Sigma_a\bg\Phi_0(x)+ \bl Q\bg,
\end{align}
with
\begin{align}\label{3.33}
\Phi(x) = \int_{-1}^1\bl\Psi\bg(x,\mu)d\mu = \int_{-1}^1[p_1\Psi_1(x,\mu)+p_2\Psi_2(x,\mu)]d\mu = \Phi_0(x) + O(\varepsilon).
\end{align}
Notice that $\Phi_0(x)$ is the leading-order estimate of the scalar flux as $\varepsilon\rightarrow 0$. Moreover, the diffusion coefficient in Eq.\ (\ref{3.32}) is identical to the one in Eq.\ (\ref{2.12}), confirming our expectation that the ALP equations yield the correct diffusion equation when $\eta$ is chosen to be of $O(1/\varepsilon)$.

In the next section we propose a choice for $\eta$ motivated by this analysis, and present numerical results that confirm the theoretical predictions and validate the ALP model.

\section{Numerical Results}\label{sec4}
\setcounter{section}{4}
\setcounter{equation}{0} 

In this section, we present numerical results that confirm the asymptotic analysis and assess the accuracy of the ALP model. To perform this assessment, we generate ensemble-averaged {\em benchmark} numerical results to the transport problem given in Eqs.\ (\ref{2.1}). We compare these results with those obtained with the ALP model, showing that it greatly outperforms the standard LP model in diffusive systems.

We also include numerical results for non-diffusive problems, in order to investigate the accuracy of the ALP model away from the diffusive limit. For these problems, we compare the benchmark results with both LP formulations {\em and} with the atomic mix model, as given 
by Eq.\ (\ref{2.13}). 

The benchmark results are attained with the benchmark method introduced in \cite{adams}, as follows. We obtain a physical realization of the system by sampling the thicknesses of individual layers from an exponential distribution with the mean values $\lambda_1$ and $\lambda_2$. This process yields $\Sigma_t(x)$, $\Sigma_s(x)$, and $Q(x)$ as histograms for this physical realization. We solve Eqs.\ (\ref{2.1}) numerically for this realization, using (i) the standard discrete ordinate method with a 16-point Gauss-Legendre quadrature set ($S_{16}$); and (ii) simple diamond differencing \cite{duderstadt} for the spatial discretization.

We repeat this process a large number of times for each problem, and average the scalar fluxes to obtain the ensemble-averaged scalar flux. Following \cite{vasques04, larsen05}, we use the Central Limit Theorem \cite{milton} to guarantee that the statistical relative error in these ensemble-averaged benchmark results is less then 1$\%$ with 95$\%$ confidence.

Numerical results for the standard LP and for the ALP models, respectively given by Eqs.\ (\ref{3.3}) and (\ref{3.5}), were also  generated using the $S_{16}$ discrete ordinate method and diamond differencing. To solve Eqs.\ (\ref{3.5}), we have defined the factor $\eta$ as
\begin{align}\label{4.1}
\eta = \left(\frac{\bl\Sigma_t\bg}{\bl\Sigma_a\bg}\right)^{1/2},
\end{align}
for $\bl\Sigma_a\bg\neq 0$. This choice of $\eta$ has the following desirable qualities:
\begin{itemize}
\item[I.] $\eta$ is $O(1/\varepsilon)$ in the 1-D diffusive systems considered in this work. This means that the ALP equations have the correct asymptotic behavior, as shown in Section \ref{sec3.2}.
\item[II.] $\eta=1$ when $\bl\Sigma_a\bg = \bl\Sigma_t\bg$. This means that the ALP equations preserve the standard LP equations in purely absorbing media.
\item[III.] $\eta$ is simple and easy to obtain. Solving the ALP equations requires the same amount of work as solving the standard LP equations. The Markov transition functions are rescaled by only one order of magnitude.
\end{itemize}

For the test problems included in this paper, we have chosen material 2 to be a void; this choice does not violate any of our previous assumptions. In fact, solid-void random mixtures correspond to known important physical applications, such as Pebble Bed Reactor cores \cite{koster,kadak} and atmospheric clouds (cf. \cite{marshak}).

\subsection{Diffusive Systems}\label{sec4.1}

To simulate the 1-D random diffusive system used in our asymptotic analysis, we consider a binary random system with total width given by
\begin{equation}
2X = (\lambda_1+\lambda_2)M = \text{total width of the system}.
\end{equation}
Here, $\lambda_1$ and $\lambda_2$ are defined in Eq.\ (\ref{3.2a}), and $M$ is given by
\begin{align}
M = \frac{1}{\varepsilon}.
\end{align}
The parameters at each spatial point $x$ in material $i$ are given by
\begin{align}
\Sigma_t(x) = \Sigma_{ti},
\,\,\,\,\,\Sigma_a(x) = \frac{\sigma_{ai}}{M^2},
\,\,\,\,\,  Q(x) = \frac{q_i}{M^2},
\end{align}
and vacuum boundary conditions are assigned at $x=\pm X$. 

\begin{table}[!ht]
\small
\centering
\caption{Parameters for Diffusive Problems}\label{tab1}
\begin{tabular}{ccc||cccccc}\hline\hline
Set & $\lambda_1$ & $\lambda_2$ & $\Sigma_{t1}$ & $\sigma_{a1}$ & $q_1$ & $\Sigma_{t2}$   & $\sigma_{a2}$  & $q_2$\\
\hline\hline
A & 1.0 & 0.5 &  &  &  & &  & \\ 
\cline{1-3}
B & 1.0 & 1.0 & 1.0 & 0.1 & 0.2 & 0 & 0 & 0\\ 
\cline{1-3}
C & 0.5 & 1.0 &   &   &   &   &   &  \\ 
\hline\hline
\end{tabular}
\end{table}
We consider the three sets of problems with parameters given in Table 1. These parameters are in agreement with the assumptions of our asymptotic analysis: $\lambda_i, \Sigma_{ti}, \sigma_{ai}$, and $q_i$ are all $O(1)$ constants, $\Sigma_{ai}$ and $Q_i$ are $O(\varepsilon^2)$, and $X$ is $O(1/\varepsilon)$.
As $\varepsilon$ decreases, the 1-D system approaches the diffusive limit. Therefore, as $M$ increases, 
we expect the solutions of both the standard and Adjusted LP models to converge to the solutions of their corresponding diffusion formulations.

Since our asymptotic analysis does not include boundary conditions, we solve the diffusion equations using the extrapolated endpoint boundary conditions
\begin{subequations}
\begin{align}
\Phi_0\left(X+\frac{2\beta}{3\bl\Sigma_t\bg}\right) = \Phi_0\left(-X-\frac{2\beta}{3\bl\Sigma_t\bg}\right)=0
\end{align}
for Eq.\ (\ref{3.16}), and
\begin{align}
\Phi_0\left(X+\frac{2}{3\bl\Sigma_t\bg}\right) = \Phi_0\left(-X-\frac{2}{3\bl\Sigma_t\bg}\right)=0
\end{align}
\end{subequations}
for Eq.\ (\ref{3.32}).

\begin{figure}
    \centering
        \includegraphics[scale=0.8]{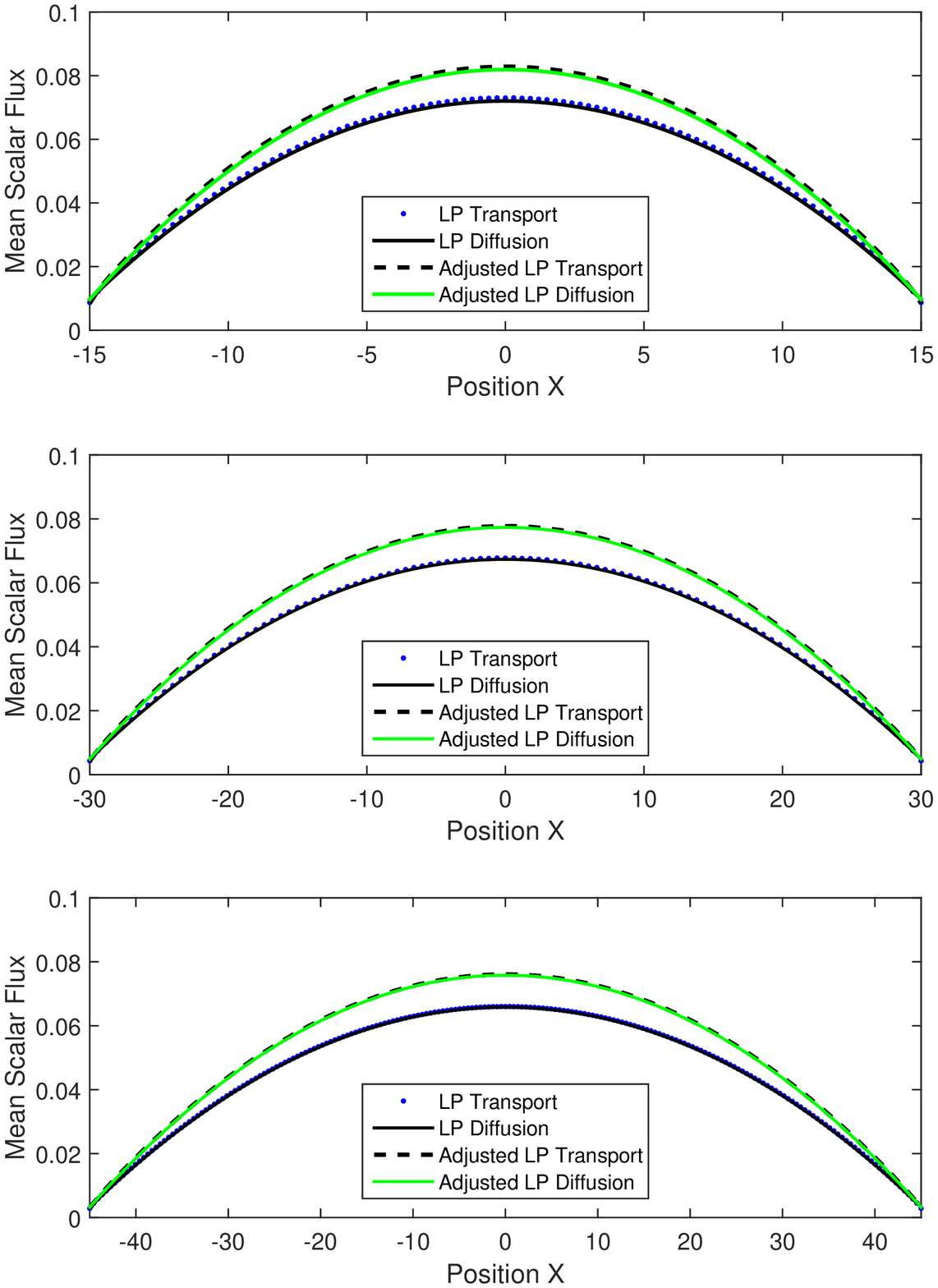}
            \caption{Solutions for the transport and diffusion formulations of the LP models for problem set A: $M = 20$ (top), $40$ (middle), $60$ (bottom)}\label{fig2}
\end{figure}

\begin{figure}
    \centering
        \includegraphics[scale=0.8]{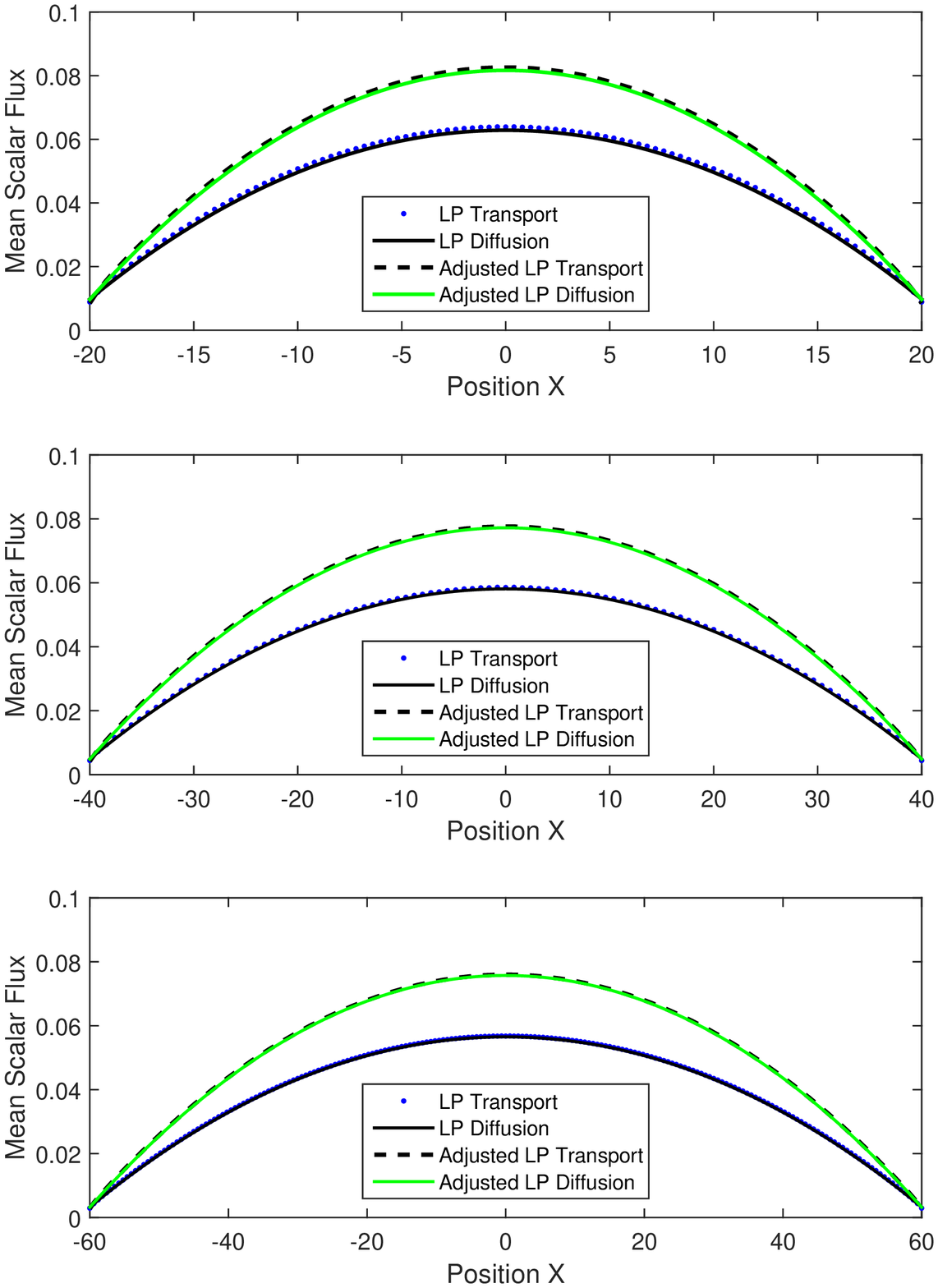}
            \caption{Solutions for the transport and diffusion formulations of the LP models for problem set B: $M = 20$ (top), $40$ (middle), $60$ (bottom)}\label{fig3}
\end{figure}

\begin{figure}
    \centering
        \includegraphics[scale=0.8]{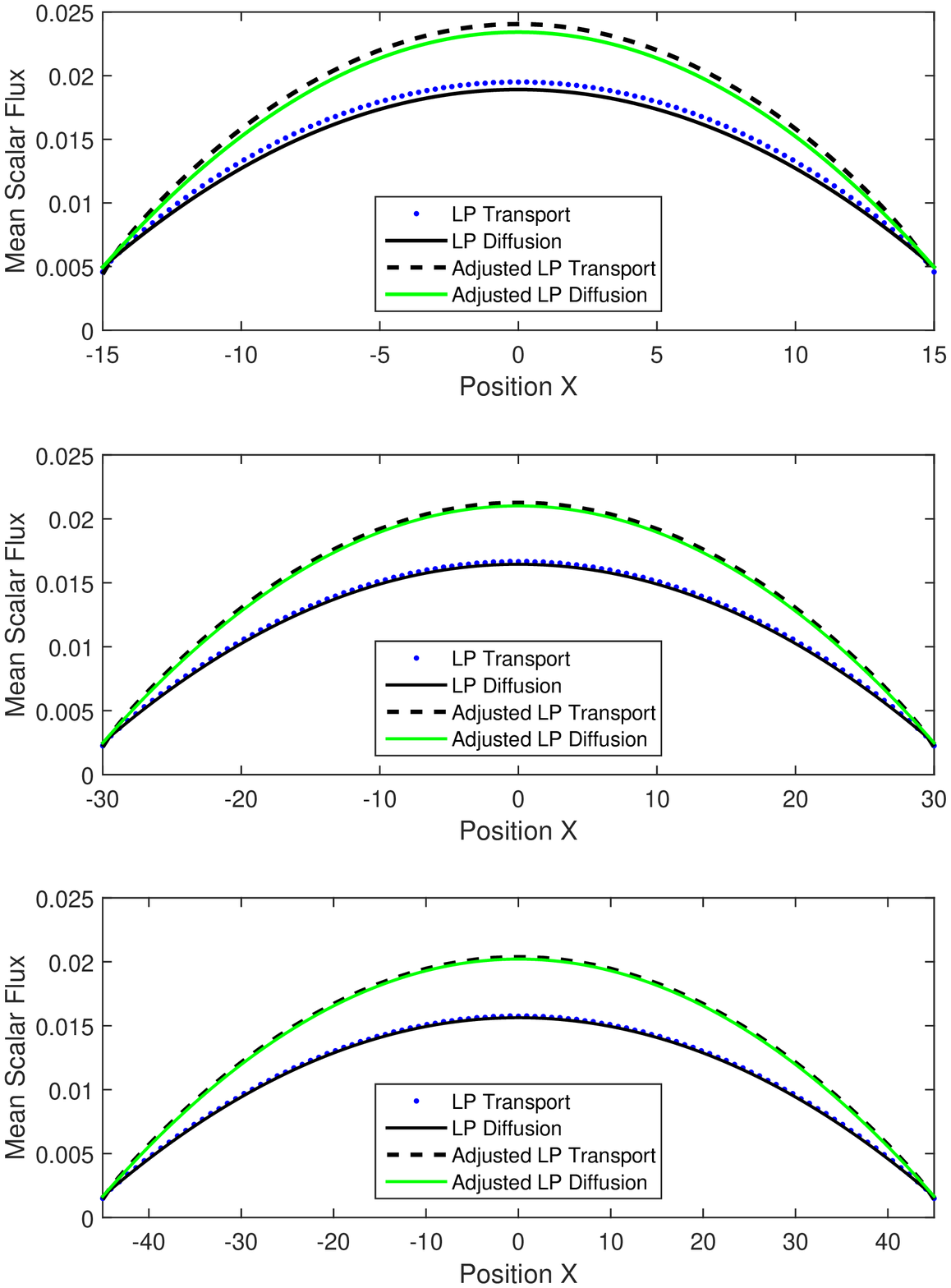}
            \caption{Solutions for the transport and diffusion formulations of the LP models for problem set C: $M = 20$ (top), $40$ (middle), $60$ (bottom)}\label{fig4}
\end{figure}

The solutions of Eqs.\ (\ref{3.3}) [LP Transport], (\ref{3.5}) [Adjusted LP Transport], (\ref{3.16}) [LP Diffusion], and (\ref{3.32}) [Adjusted LP Diffusion] are plotted in Figures \ref{fig2}, \ref{fig3}, and \ref{fig4}. The expected convergence as $M$ increases is clearly seen in each figure, confirming the results of our asymptotic analysis.

\begin{figure}
    \centering
        \includegraphics[scale=0.8]{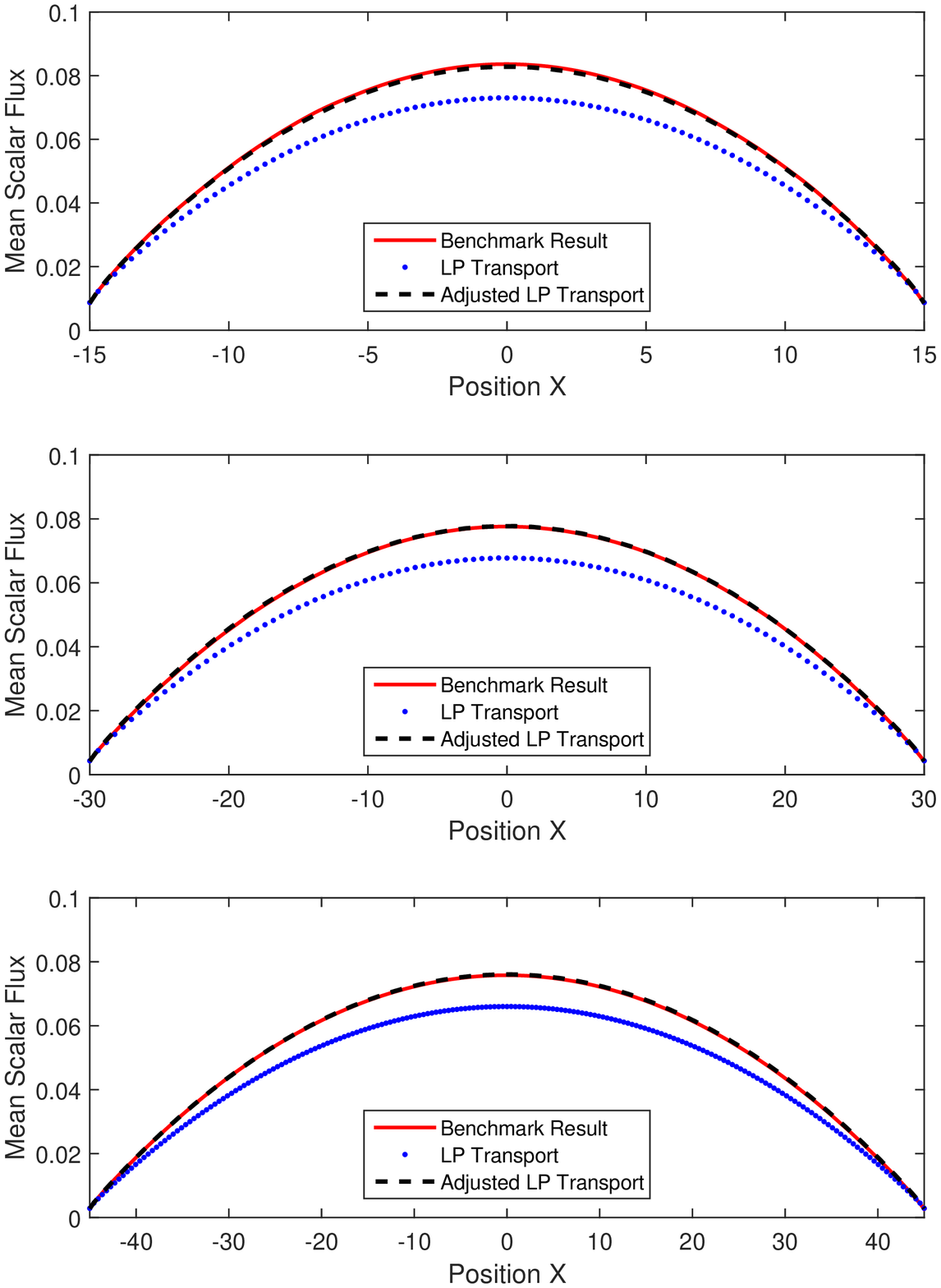}
            \caption{Comparison between LP and benchmark results for problem set A: $M = 20$ (top), $40$ (middle), $60$ (bottom)}\label{fig5}
\end{figure}

\begin{figure}
    \centering
        \includegraphics[scale=0.8]{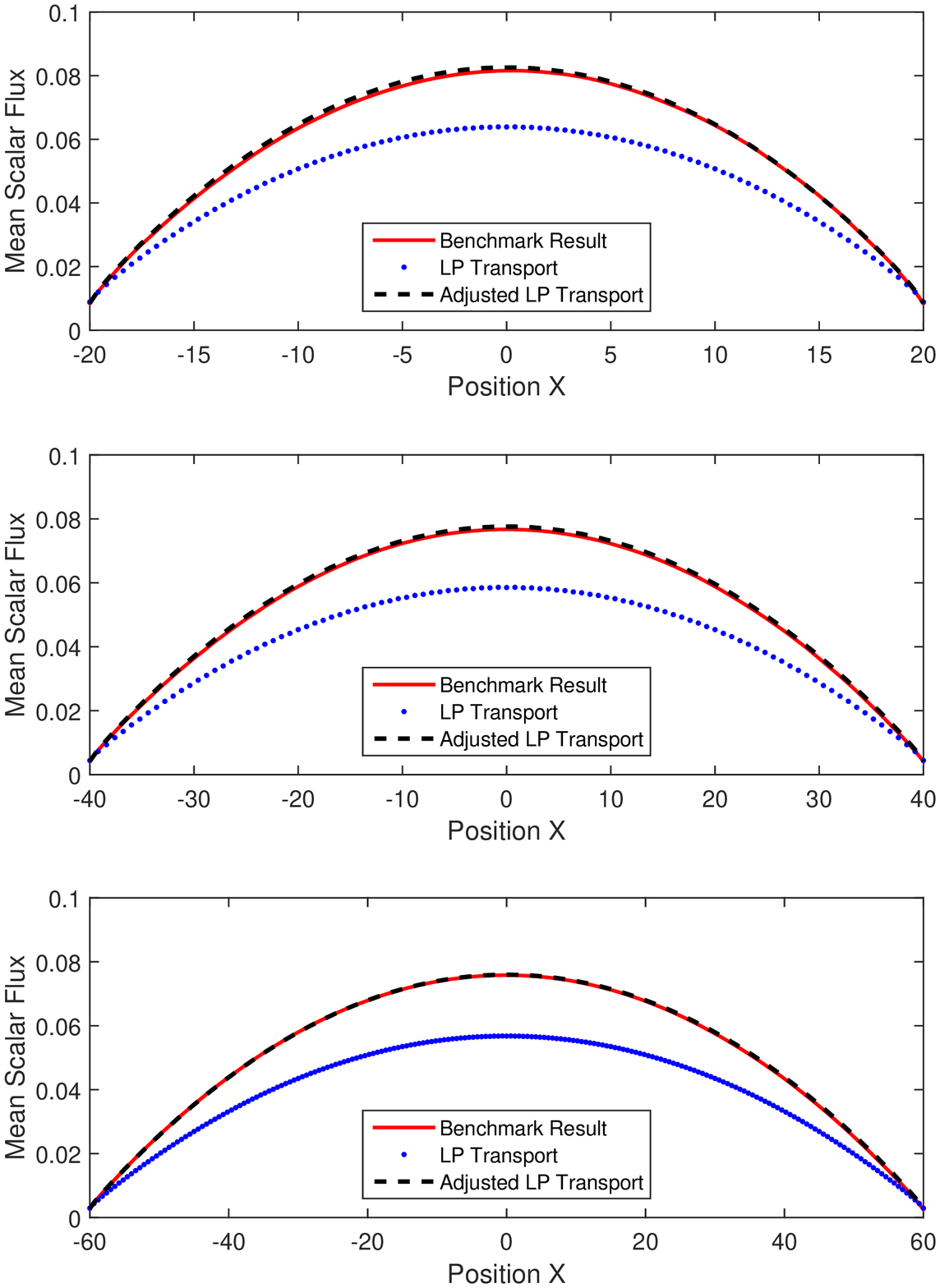}
            \caption{Comparison between LP and benchmark results for problem set B: $M = 20$ (top), $40$ (middle), $60$ (bottom)}\label{fig6}
\end{figure}

\begin{figure}
    \centering
        \includegraphics[scale=0.8]{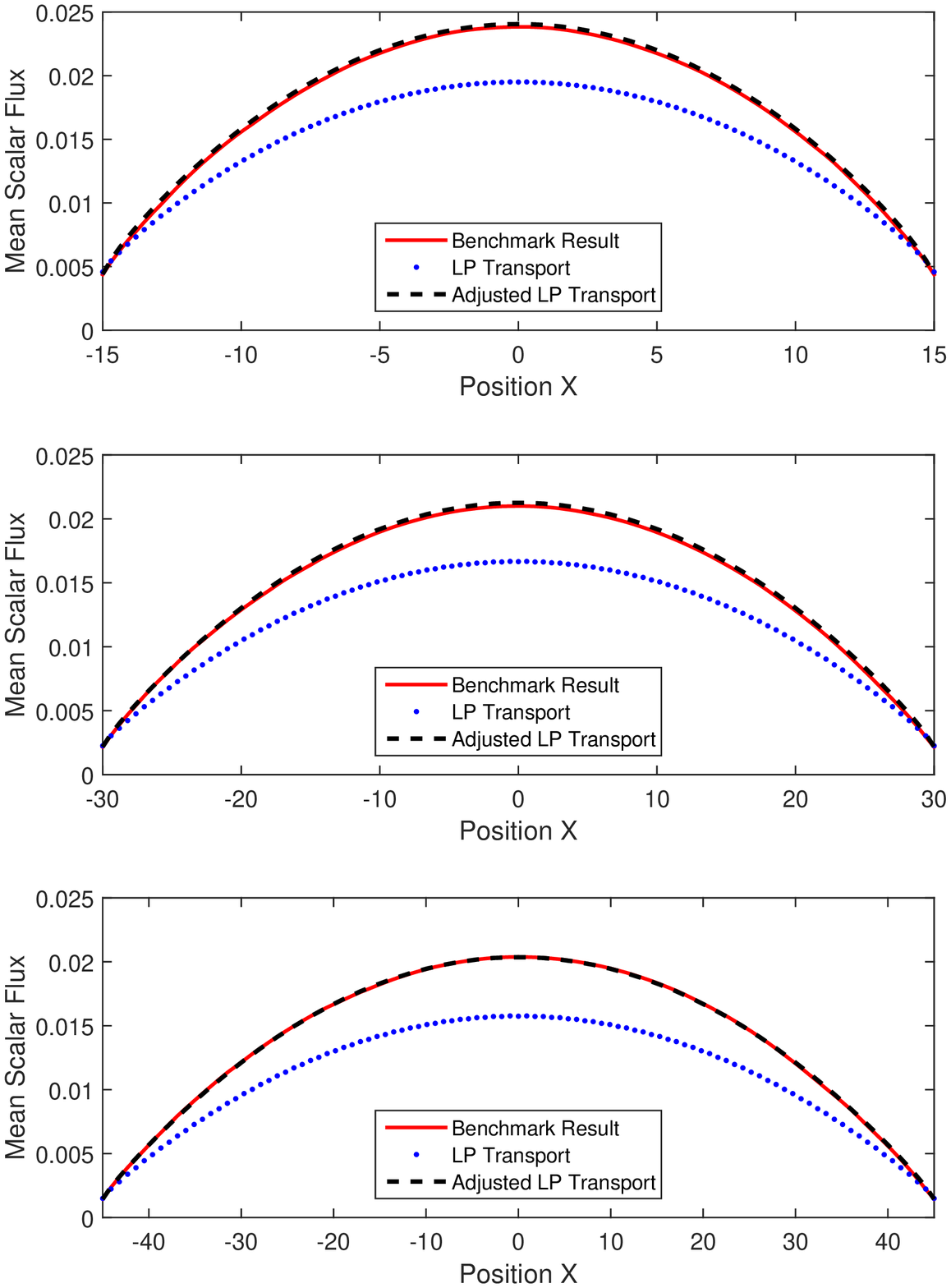}
            \caption{Comparison between LP and benchmark results for problem set C: $M = 20$ (top), $40$ (middle), $60$ (bottom)}\label{fig7}
\end{figure}

In Figures \ref{fig5}, \ref{fig6}, and \ref{fig7}, we again plot the solutions of Eqs.\ (\ref{3.3}) and (\ref{3.5}), this time comparing them with the benchmark numerical results of Eqs.\ (\ref{2.1}). The standard LP model systematically disagrees with the benchmark results. As predicted, its estimate for the scalar flux is incorrectly flattened. On the other hand, the ALP model estimates the correct scalar flux with great accuracy.

\begin{figure}
    \centering
        \includegraphics[scale=0.8]{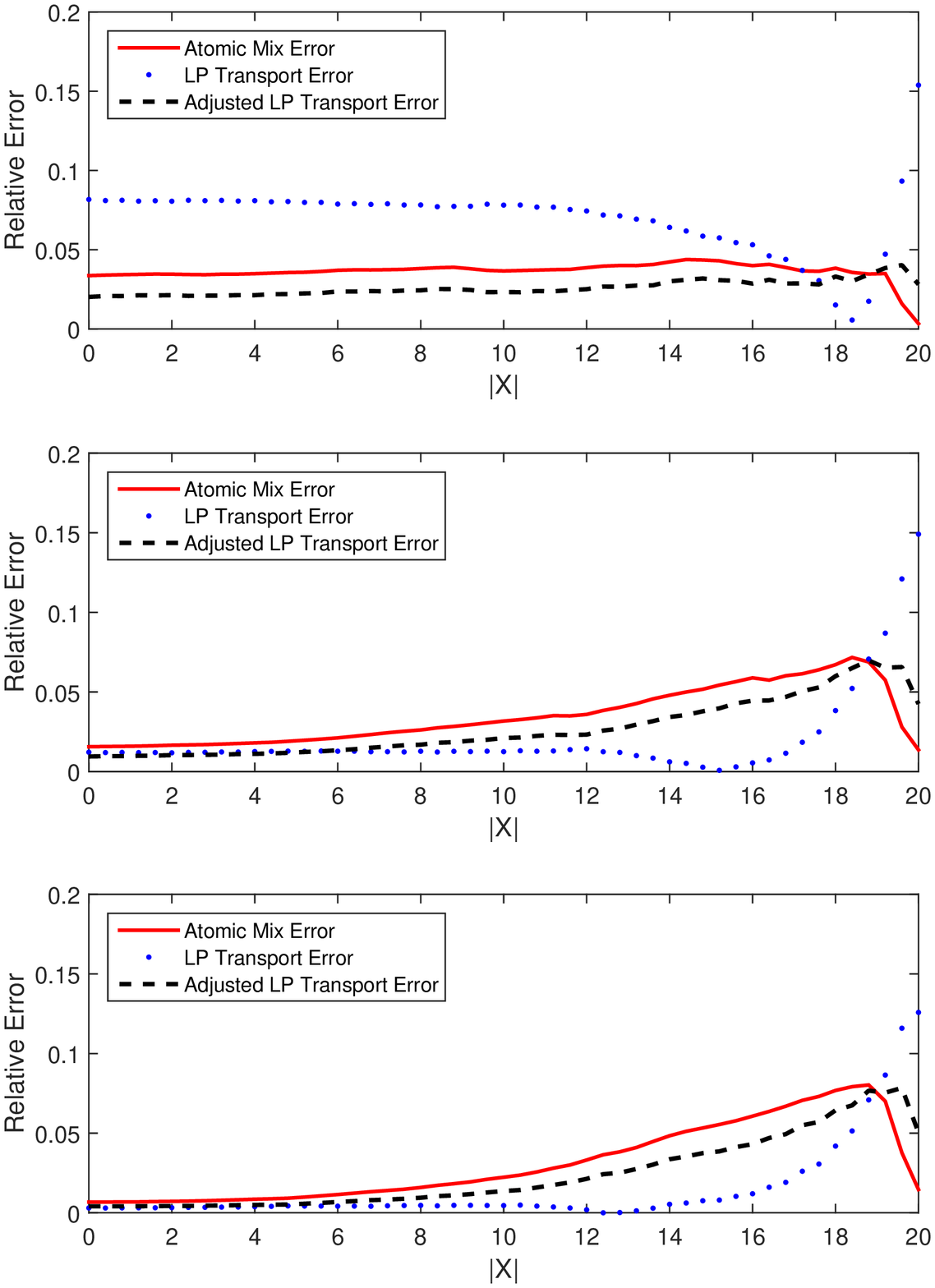}
            \caption{Absolute values of the LP and atomic mix relative errors with respect to the benchmark solutions for problem set D: $\Sigma_{s1} = 0.99$ (top), $0.95$ (middle), $0.9$ (bottom)}\label{fig8}
\end{figure}

\begin{figure}
    \centering
        \includegraphics[scale=0.8]{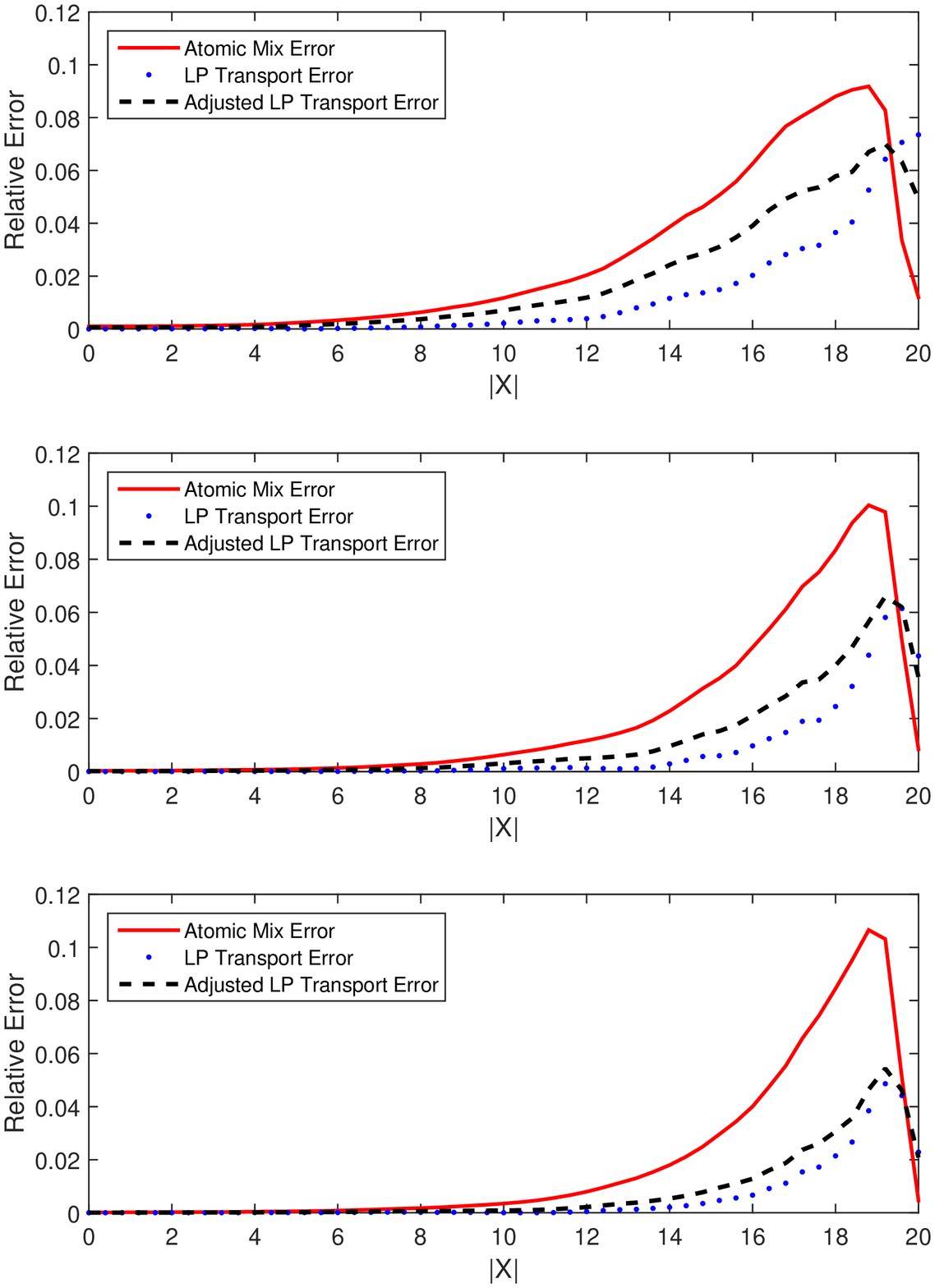}
            \caption{Absolute values of the LP and atomic mix relative errors with respect to the benchmark solutions for problem set E: $\Sigma_{s1} = 0.7$ (top), $0.5$ (middle), $0.3$ (bottom)}\label{fig9}
\end{figure}

\begin{figure}
    \centering
        \includegraphics[scale=0.8]{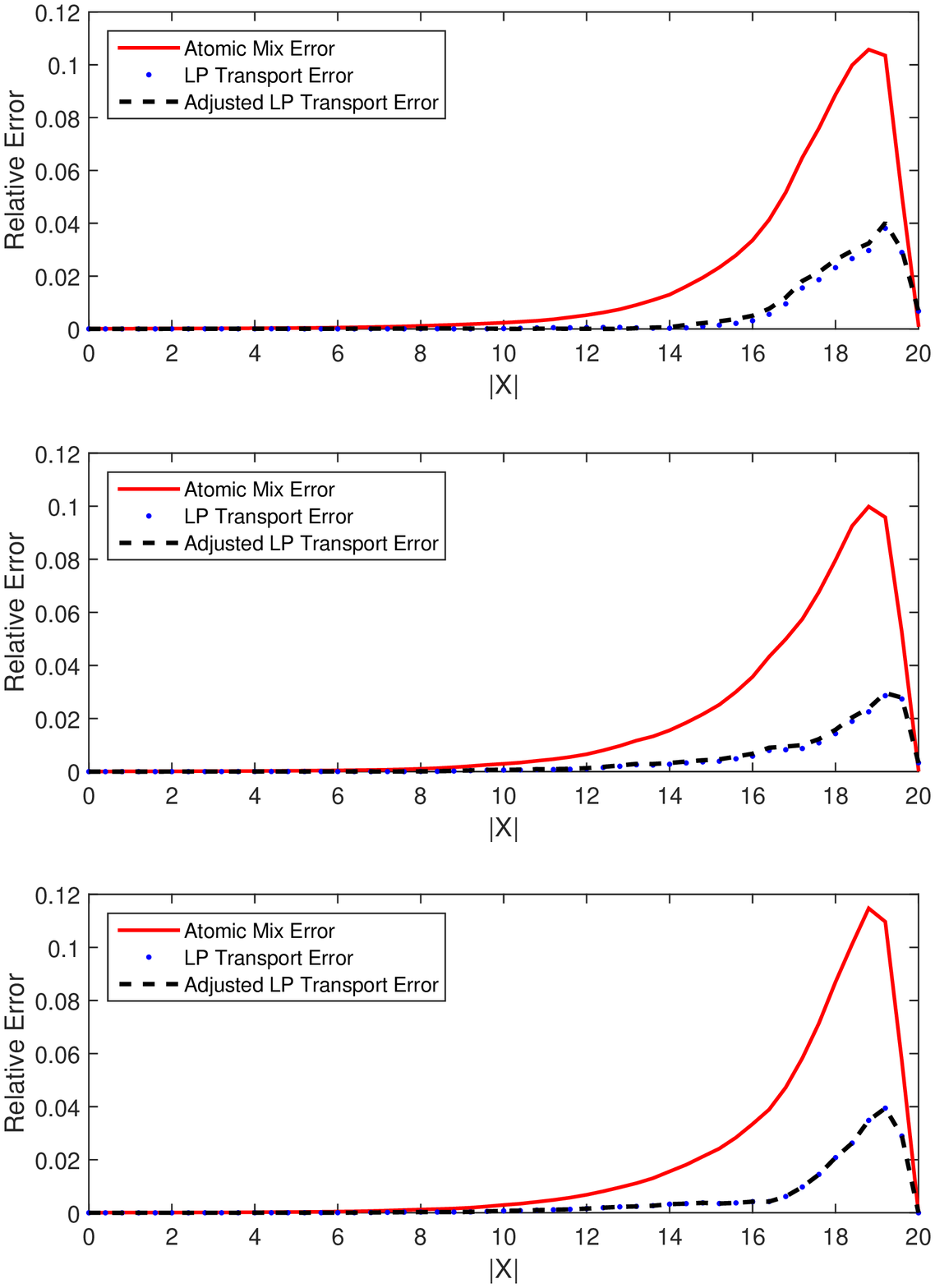}
            \caption{Absolute values of the LP and atomic mix relative errors with respect to the benchmark solutions for problem set F: $\Sigma_{s1} = 0.1$ (top), $0.05$ (middle), $0.0$ (bottom)}\label{fig10}
\end{figure}

For a better analysis of these results, we define the relative errors of the models with respect to the benchmark solutions as
\begin{subequations}\label{4.6}
\begin{align}
Err^{(LP)}(x) &= \frac{\Phi^{(LP)}(x)-\Phi^{(B)}(x)}{\Phi^{(B)}(x)},\label{4.6a}\\
Err^{(ALP)}(x) &= \frac{\Phi^{(ALP)}(x)-\Phi^{(B)}(x)}{\Phi^{(B)}(x)}.\label{4.6b}
\end{align}
Here, $\Phi^{(B)}(x)$, $\Phi^{(LP)}(x)$, and $\Phi^{(ALP)}(x)$ are the estimated scalar fluxes obtained with the benchmark, standard LP, and Adjusted LP models, respectively. The numerical estimates of each model for the scalar flux at $x=0$ are given in Table \ref{tab2}, as well as the correspondent (percent) relative errors.  
\begin{table}[!ht]
\small
\centering
\caption{Ensemble-averaged scalar fluxes and relative errors at x = 0 (Diffusive Problems)}\label{tab2}
\begin{tabular}{c|c|ccc|cc}\hline\hline
Set & $M$ & $\Phi^{(B)}$ & $\Phi^{(LP)}$ & $\Phi^{(ALP)}$ & $Err^{(LP)}$ & $Err^{(ALP)}$   \\
\hline\hline
& 20 & 0.0836 & 0.0730 & 0.0828 & -12.68\% & -0.96\%\\
 \cline{2-7}
 A & 40 & 0.0776 & 0.0677 & 0.0777 & -12.76\% & 0.13\%\\
 \cline{2-7}
 &60& 0.0758 & 0.0660 & 0.0759 & -12.93\% & 0.13\%\\ 
\hline
 & 20 & 0.0816 & 0.0639 & 0.0825 & -21.69\% & 1.10\%\\
 \cline{2-7}
 B & 40 & 0.0767 & 0.0585 & 0.0776 & -23.73\% & 1.17\% \\
 \cline{2-7}
 & 60 & 0.0758 & 0.0567 & 0.0759 & -25.20\% & 0.13\%\\
 \hline
  & 20 & 0.0238 & 0.0195 & 0.0239 & -18.07\% & 0.42\%\\
 \cline{2-7}
 C& 40 &0.0210 & 0.0167 & 0.0213 & -20.48\% & 1.43\%\\
  \cline{2-7}
&  60 & 0.0204 & 0.0157 & 0.0204 & -23.04\% & $\approx 0$\\
 \hline\hline
\end{tabular}
\end{table}

As anticipated, the solutions plotted in Figures \ref{fig5}, \ref{fig6}, and \ref{fig7} and the data given in Table \ref{tab2} confirm our claim: the ALP model greatly outperforms its standard counterpart in the diffusive limit. Next, we investigate the performance of the ALP model in non-diffusive systems.

\subsection{Non-Diffusive Systems}\label{sec4.2}

Here, we present numerical results to assess the accuracy of the ALP model away from the diffusive limit discussed earlier. As in Section \ref{sec4.1}, we compare the LP and ALP solutions with benchmark numerical solutions of the transport problem given by Eqs.\ (\ref{2.1}). Moreover, we provide numerical results for the {\em atomic mix approximation} given by Eq.\ (\ref{2.13}).

We define the 1-D binary random systems such that:
\begin{itemize}
\item[\textbf{I}] The total width of the system is $40$;
\item[\textbf{II}] The mean width of layers is given by $\lambda_1=\lambda_2 =1.0$;
\item[\textbf{III}] Vacuum boundary conditions are assigned at $x=\pm 20$;
\item[\textbf{IV}] Material 2 is a void: $\Sigma_{t2}=\Sigma_{s2}=Q_2=0$;
\item[\textbf{V}] Material 1 has total cross section $\Sigma_{t1}=1.0$ and source $Q_1=0.2$;
\item[\textbf{VI}] The scattering cross section $\Sigma_{s1}$ is given in Table \ref{tab3} for three different sets of problems.
\end{itemize}
\begin{table}[!ht]
\small
\centering
\caption{Scattering Cross Sections for Non-Diffusive Problems}\label{tab3}
\begin{tabular}{cccc}\hline\hline
 &\multicolumn{3}{c}{Scattering cross section $\Sigma_{s1}$} \\ \cline{2-4}
Set & Choice 1 & Choice 2 & Choice 3\\
\hline\hline
D & 0.99& 0.95& 0.9\\ 
E & 0.7& 0.5& 0.3\\ 
F & 0.1& 0.05& 0.0\\ 
\hline\hline
\end{tabular}
\end{table}

Similarly to Eqs.\ (\ref{4.6a}) and (\ref{4.6b}), we define the relative error of the atomic mix model with respect to the benchmark solutions as
\begin{align}
Err^{(AM)}(x) &= \frac{\Phi^{(AM)}(x)-\Phi^{(B)}(x)}{\Phi^{(B)}(x)},
\end{align}
\end{subequations}
where $\Phi^{(AM)}(x)$ is the estimated scalar flux obtained by solving Eq.\ (\ref{2.13}). The {\em absolute values} of the relative errors given by Eqs.\ (\ref{4.6}) for the problems sets D, E, and F are given in Figures \ref{fig8}, \ref{fig9}, and \ref{fig10}, respectively. In Table \ref{tab4}, we provide the numerical values obtained with each model for the scalar flux at $x=0$, as well as the correspondent (percent) relative errors.

The $x$-axis in the plots shown in Figures \ref{fig8}-\ref{fig10} is the distance from the origin, due to the symmetry of the solutions. The spike in the errors close to $|x|=20$ is simply a boundary effect: as the system becomes more absorbing, the scalar flux away from the boundaries approaches the ``infinite-medium" solution $\bl Q\bg/(2\bl\Sigma_a\bg)$. 

\begin{table}[!ht]
\small
\centering
\caption{Ensemble-averaged scalar fluxes and relative errors at x = 0 (Non-Diffusive Problems)}\label{tab4}
\begin{tabular}{c|c|cccc|ccc}\hline\hline
Set & $\Sigma_{s1}$ & $\Phi^{(B)}$ & $\Phi^{(LP)}$ & $\Phi^{(ALP)}$ & $\Phi^{(AM)}$ & $Err^{(LP)}$ & $Err^{(ALP)}$ & $Err^{(AM)}$   \\
\hline\hline
& 0.99 & 13.396 & 12.302 & 13.667 & 13.847 &-8.17\% & 2.02\% &3.37\% \\
\cline{2-9}
D & 0.95 & 3.8050 & 3.7587 & 3.8412 & 3.8644 &-1.22\% &0.95\% &1.56\%\\
\cline{2-9}
 & 0.9 & 1.9729 & 1.9669 & 1.9809 & 1.9863 &-0.30\%&0.41\%&0.68\%\\
 \hline
& 0.7 & 0.6659 & 0.6659 & 0.6663 & 0.6665 &$\approx 0$&0.06\%&0.09\%\\
\cline{2-9}
E & 0.5 & 0.3999 & 0.3999 & 0.3999 & 0.4000&$\approx 0$&$\approx 0$& 0.03\%\\
\cline{2-9}
& 0.3 & 0.2857 & 0.2857 & 0.2857 & 0.2857 & $\approx 0$ & $\approx 0$ & $\approx 0$\\
\hline
&0.1 & 0.2222 & 0.2222 & 0.2222 &0.2222 & $\approx 0$ & $\approx 0$ & $\approx 0$\\
\cline{2-9}
F & 0.05 & 0.2105 & 0.2105& 0.2105& 0.2105 & $\approx 0$ & $\approx 0$ & $\approx 0$\\
\cline{2-9}
 & 0.0 & 0.2000 & 0.2000 & 0.2000 & 0.2000 & $\approx 0$ & $\approx 0$ & $\approx 0$\\
 \hline\hline
\end{tabular}
\end{table}

The first plot in Figure \ref{fig8} shows that the ALP model consistently outperforms both the standard LP and the atomic mix models for the scattering ratio $0.99$. However, as the scattering ratio decreases, standard LP overtakes the adjusted model in terms of accuracy. 

Furthermore, we see from the plots in Figure \ref{fig10} and from the data in Table \ref{tab4} that the ALP equations preserve the accuracy of the standard model for systems in which absorption is high. This validates our choice of $\eta$ in Eq.\ (\ref{4.1}). 

Finally, it is clear that the proposed choice of $\eta$ does not improve on the LP results for the whole spectrum of the scattering ratio. Nevertheless, it significantly improves the LP equations for diffusive systems, and even outperforms the atomic mix model when the diffusive parameters are slightly relaxed. This result paves the road to finding an expression for $\eta$ that works for the whole spectrum of the scattering ratio, which would be of extreme importance for several random media applications. 

\section{Discussion}\label{sec5}
\setcounter{section}{5}
\setcounter{equation}{0} 

This paper presents an adjustment to the standard Levermore-Pomraning equations for diffusive problems in slab geometry.    
This adjustment is motivated by an analysis that shows that, under certain diffusive conditions, the asymptotic behavior of the LP model deviates from the correct result. This is confirmed by numerical simulations, providing an explanation for previously-observed inaccuracies in the LP model for diffusive problems.

The asymptotic analysis is valid for physical systems that (i) have weak absorption and sources; (ii) are optically thick; and (iii) consist of a large number of
material layers with mean thicknesses comparable to (or small compared to) a mean free path. It is not
necessary that the system be highly-scattering at all points; void regions are permitted.

The analysis shows that, by introducing a factor $\eta = O(1/\varepsilon)$ to the original LP model, one can fix its asymptotic behavior. In this work, we have chosen this factor to be $\eta = \left(\bl\Sigma_t\bg/\bl\Sigma_a\bg\right)^{1/2}$, which (i) yields very accurate results in the diffusive limit; (ii) preserves the exactness of the LP model for purely absorbing problems; (iii) is simple to compute. 

Claims (i) and (ii) above are confirmed by our numerical results. The numerical simulations presented in this paper take place in solid-void random mixtures. These correspond to important physical applications, such as neutron diffusion in Pebble Bed Reactor (PBR) cores \cite{koster,kadak}, and radiative transfer in atmospheric clouds (cf. \cite{marshak}). In PBR cores, the random structure consists of a mixture of fuel pebbles (solid) and Helium (void). In atmospheric clouds, the system is a random mixture of water droplets (solid) and air (void).

From a qualitative viewpoint, the factor $\eta$ introduced in the adjusted model consists of a rescaling of the Markov transition functions: by substituting the original $O(1)$ mean widths $\lambda_i$ for the ``rescaled" widths $\lambda_i/\eta = O(\varepsilon)$, we are solving the LP equations in an artificial atomic mix limit, which yields the asymptotically correct diffusion formulation given by Eq.\ (\ref{2.12}).

Asymptotic diffusion approximations of the LP equations have been considered before \cite{pomraning91}. However, these were performed for systems in which material layers are optically thick - that asymptotic limit is fundamentally different from the
one considered here. Moreover, to our knowledge, corresponding asymptotic limits have not been applied to the original transport equation, so it is unclear whether the results are physically correct. 

Although the analysis developed in this paper can be extended to multidimensional problems, the diffusion equations obtained no longer generally hold. For instance, in general 3-D problems, spatial correlations to the cross sections will lead to anisotropic diffusion. Nevertheless, for the diffusive systems considered here, the ALP equations are a clear improvement over their standard counterpart, and provide a theoretical tool for examining future generalizations of the LP model. The numerical results
indicate that the adjusted model remains a valid alternative to current methods even when the diffusive parameters are slightly relaxed. The correction proposed in this paper can be taken as a first step in achieving a generalization of the LP equations that is accurate for the whole spectrum of the scattering ratio. We intend to refine this idea in future work.

\section*{Acknowledgments}
This work was supported by funds from the German Federal Ministry of Education and Research under grant 02S9022A. The responsibility for the content of this publication lies with the authors.


\pagebreak


\pagebreak

\end{document}